\documentstyle{article}
\def\P{{\rm I\hskip -2pt P}}
\def\N{{\rm I\hskip -2pt N}}
\def\R{{\rm I\hskip -2pt R}} 
\def\Aff{{\rm Aff}}
\def\tenreg{\rm\bf}
\def\findemo{{\hfill QED}\vskip 0.9cm}
\newtheorem{thm}{\bf Theorem}[section]
\newtheorem{cor}[thm]{\bf Corollary}
\newtheorem{lem}[thm]{\bf Lemma}
\newtheorem{prop}[thm]{\bf Proposition}

\newtheorem{defn}[thm]{\bf Definition}

\def\vtab{\vspace{0.3cm}}
\def\dd#1{{\partial\over\partial #1}}
\def\hfl#1#2{\smash{\mathop{\hbox to 8mm{\rightarrowfill}}
\limits^{\scriptstyle #1}_{\scriptstyle #2}}}
\def\hfll#1#2{\smash{\mathop{\hbox to 8mm{\leftarrowfill}}
\limits^{\scriptstyle #1}_{\scriptstyle #2}}}
\def\vfl#1#2{\llap{$\scriptstyle #1$}\left\downarrow
\vbox to 4mm{}\right.\rlap{$\scriptstyle #2$}}

\begin{document}

 \title{Duality of subanalytic sets}
 
 \author{Fran\c{c}ois Pointet\\
         Institut de math\'ematiques\\
            Universit\'e de Lausanne\\
            CH-1015 Lausanne-Dorigny\\
            e-mail: francois.pointet@ima.unil.ch}
 
 \date{}
 \maketitle

\begin{abstract}
We study the link between a compact hypersurface
in $\P^{n+1}$ and the set of all its tangent planes. In this context, we identify
$\P^{n+1}$ to the set of linear subspaces of codimension one by orthogonal
complementarity. This gives rise to a kind of duality which has already been studied in
\cite{Bru1}, \cite{Rom} and relates a hypersurface to the set of its tangent planes. But in
these papers the dual, in this sense, of the set of tangent planes of a hypersurface was
not defined and iteration of the procedure was not possible. Therefore we extend this
type of duality to more general sets and achieve a procedure which can be iterated and
gives in fact an involution.

\end{abstract}

\section{Introduction}
The purpose of this work is to generalize a notion of duality, which is well known in the 
case of algebraic geometry \cite{GKZ}, to the analytic case. 

Let us recall what is duality in the algebraic case.
We denote by $\P^n$ the $n$-dimensional real projective space and by
$\pi\,:\R^{n+1}\setminus\{0\}\rightarrow\P^n$ the canonical projection. A $(n-1)$-dimensional
linear subspace $L\subset\P^{n}$ (or hyperplane of $\P^n$) is by definition the image through
$\pi$ of a vector subspace $V\subset\R^{n+1}$ of codimension 1. We denote by
$(\P^n)^*$ the space of codimension 1 linear subspaces of $\P^n$. There is an isomorphism 
$\varphi$ from $(\P^n)^*$ to $\P^n$ defined by $\varphi(L)=\pi(V^{\perp}\setminus\{0\})$,
where $V$ is such that $\pi(V\setminus\{0\})=L$ and
$V^{\perp}$ is the normal vector space to $V$ with respect to the usual scalar product in 
$\R^{n+1}$. Now we can consider
$M\subset\P^n$ an algebraic hypersurface (the zero set of an homogeneous polynomial) and
$M^*\subset(\P^n)^*$ the set of its tangent hyperplanes. Then $M^*$, which is called the dual
of $M$, is an algebraic hypersurface of
$(\P^n)^*$ and $(M^*)^*=M$ via the canonical isomorphism from $(\P^n)^{**}$ to $\P^n$. 
We can see \cite{GKZ} to have some results we can obtain in the algebraic case.

The motivation in generalizing this study to non-algebraic hypersurfaces appears when we try 
to reconstruct hypersurfaces in $\R^{n+1}$ using outlines \cite{GW}, \cite{CB} and
\cite{Poi}. First we give our approach to study the dual of a hypersurface in $\R^{n+1}$ and
we give the results we obtain by generalizing the notion of duality.

We denote by $\Aff(n,1)$ the set of affine hyperplanes in $\R^{n+1}$.
The affine dual of a smooth hypersurface $M$ in $\R^{n+1}$ is the set of its affine tangent 
planes \cite{Bru1}, \cite{Rom}. As we want to obtain an iterable process, it is
important to see the dual of $M$ as a subset of $\P^{n+1}$. Let us denote by
$i\,:\R^{n+1}\rightarrow\P^{n+1}$ the injection given by $i(x)=[x,1]$. If $\pi\in\Aff(n,1)$,
then we have that $\overline{i(\pi)}$ is a hyperplane in $\P^{n+1}$. So we can see
$\Aff(n,1)$ as a part of
$(\P^{n+1})^*$, in fact $\Aff(n,1)$ is isomorphic to $(\P^{n+1})^*\setminus\{\sigma\}$ where 
$\sigma$ is the hyperplane given by $\{[x,0]\mid x\in\R^{n+1}\setminus\{0\}\}$. Now using the
map $\varphi$ given above, we can identify $\Aff(n,1)$ with
$\P^{n+1}\setminus\{[0_{n+1},1]\}$, and we have an expression of this identification with
the map\,:
$$\matrix{\varphi&\,:&\Aff(n,1)&\longrightarrow&\P^{n+1}\setminus\{[0_{n+1},1]\}\cr
                 &&\pi&\longmapsto&[v_1(\pi),\dots,v_{n+2}(\pi)]\cr}$$
where $\pi$ is given by $v_1(\pi)x_1+\dots+v_{n+1}(\pi)x_{n+1}=v_{n+2}(\pi)$.{\it }
Using this identification, we can define the dual of a hypersurface in $\R^{n+1}$.
\begin{defn}
\rm
Let $M$ be a hypersurface of $\R^{n+1}$, we define the {\it affine Gauss map of $M$} 
to be\,: 
$$\matrix{G^{\Aff}_M&\,:&M&\longrightarrow&\P^{n+1}\cr &&m&\longmapsto&[N(m),<N(m),m>]\cr}$$
where $N(m)$ is a normal to $M$ at $m$ and $<,>$ denotes the usual scalar product in $\R^{n+1}$. The set
$G^{\Aff}_M(M)$ is called the {\it dual of $M$}.

We denote by {\tenreg H} the set of smooth hypersurfaces in $\R^{n+1}$ and we define the map\,:
$$\matrix{G^{\Aff}&\,:&\hbox{\tenreg H}&\longrightarrow&\hbox{\tenreg
P}\,(\P^{n+1})\cr &&M&\longmapsto&G^{\Aff}_M(M)\cr}$$ where $\hbox{\tenreg P}\,(\P^{n+1})$
denotes the power set of $\P^{n+1}$.
\end{defn}
The map $G^{\Aff}$ associates to every hypersurface its dual, we can see \cite{Bru1} to have a study of
the singularities we can find in $G^{\Aff}_M$ for a generic hypersurface $M$.

The problem in this point of view comes from the fact that the process of taking
the dual is not iterable. So in this work, we are going to generalize the notion of duality to a good
family of subsets in $\P^{n+1}$ to have an iterable process. We will obtain the following theorem
(where {\tenreg Ad} denotes the set of admissible subsets of $\P^{n+1}$, a notion that will be defined
later)\,:

\vtab
{\bf Theorem \ref{bigthm}.} The map $G^{\Aff}$ restricted to {\tenreg Ad} is an involution.

\vtab
After this generalization, we can come back to the motivating problem. Let us recall what is an outline
and what we mean by reconstruction of hypersurfaces using outlines.

\begin{defn}
\rm
Let $M$ be a hypersurface of $\R^{n+1}$ and $\sigma\in\P^{n+1}$, we denote by
$\pi_{\sigma}\,:\R^{n+1}\rightarrow\R^{n+1}$ the orthogonal projection of $\R^{n+1}$ on $\sigma^{\perp}$.
The {\it outline of $M$ in the direction $\sigma$} is the set
$C^M_\sigma=\pi_\sigma(\Sigma(\pi_\sigma\!\mid_M))$ where $\Sigma(\pi_\sigma\!\mid_M)$ denotes the
critical set of the restriction of $\pi_\sigma$ to $M$.
\end{defn}

The question of reconstruction is the following\,: Let $M\subset\R^{n+1}$ be a hypersurface and suppose you
have a certain number of its outlines, can we reconstruct exactly $M$\,?

In \cite{Poi} we give an answer for compact hypersurfaces supposing we have the outlines for a set of
directions with a certain property. During this study we use the notion of duality of hypersurfaces,
without studying this notion in details. Using the generalization obtained, it is easy to prove that the
dual of an outline (in the hyperplane containing the outline) of a hypersurface $M$ is contained in the
intersection of the dual of $M$ with a hyperplane of $\P^{n+1}$ \cite{BG}. So for $M$ a compact analytic
hypersurface in $\R^{n+1}$ and
$[v]\in\P^n$ a direction we obtain\,:
$$\matrix{M_{}&\hfl{G^{\Aff}\circ
i}{}&G^{\Aff}(i(M))\cr\vfl{\pi_{[v]}}{}&&\vfl{}{\cap\P([v,0]^\perp)}\cr
C^M_{[v]}&\hfll{\xi\circ G^{\Aff}}{}&G^{\Aff}(i(M))\cap\P([v,0]^\perp)\cr}$$ This gives a
new method to reconstruct hypersurfaces using outlines.

\vtab
Let us give a description of each chapter. In chapter 2 we make the first step of the generalization by
defining the dual of a subanalytic set of $\P^{n+1}$. Then we restrict the study to admissible sets to obtain an iterable process, this point is made in the chapter 3. To have
a link between duality and outlines, we study, in chapter 4, intersections of admissible sets with
hyperplanes of $\P^{n+1}$. Chapter 5 contains the elaboration of this link, and finally in chapter 6
we find an application of this duality.

This work is a part of the PhD thesis of the author \cite{Poi2}.

\section{Generalization of the affine dual}

In this chapter, we will generalize the map $G^{\Aff}$ to subanalytic subsets of $\P^{n+1}$ and we will
prove that the image of a subanalytic set through $G^{\Aff}$ is subanalytic. To do that, we will use
classical results on subanalytic sets \cite{BM}, \cite{Loj2}.

Let us denote by
$\pi\,:\R^{n+2}\setminus\{0\}\rightarrow\P^{n+1}$ the canonical projection.

\begin{defn}
\rm
We say that $L\subset\P^{n+1}$ is a {\it $k$-linear subspace of $\P^{n+1}$} if there exists
$V\subset\R^{n+2}$ a $(k+1)$-dimensional vector subspace with $L=\pi(V\setminus\{0\})$.

Let $L\subset\P^{n+1}$ be a $k$-linear subspace of $\P^{n+1}$ with $\pi(V\setminus\{0\})=L$. The {\it 
orthogonal space of $L$} is the $(n+1-k)$-linear subspace $L^{\perp}$ given by
$L^{\perp}=\pi(V^{\perp}\setminus\{0\})$, where $V^{\perp}$ is the orthogonal vector space of
$V$ with respect to the usual scalar product in $\R^{n+2}$.
\end{defn}
First we have to generalize $G^{\Aff}$ to smooth hypersurfaces of $\P^{n+1}$.
\begin{defn}
\rm
Let $M\subset\P^{n+1}$, the set $\pi^{-1}(M)$, denoted by $C_M$, is called the {\it cone of $M$}.

Let $M\subset\P^{n+1}$ be a hypersurface and $m\in M$. The set $L_mM\,:=\pi(T_{m'}C_M)$, where
$m'$ is such that $\pi(m')=m$, is called the {\it linear tangent space of $M$ at $m$}.
\end{defn}
We remark that $L_mM$ is independent of the choice of $m'$. 
\begin{defn}
\rm
Let $M\subset\P^{n+1}$ be a hypersurface, we define the {\it affine Gauss map of $M$} to be\,:
$$\matrix{G^{\Aff}_M&\,:&M&\longrightarrow&\P^{n+1}\cr &&m&\longmapsto&(L_mM)^{\perp}\cr}$$
\end{defn}
Now we can generalize the map $G^{\Aff}$ to subanalytic sets, but first we recall some
interesting properties of such sets.

\begin{prop}{\bf (\cite{BM})}

- A subanalytic set admits a Whitney stratification.

- The closure of a subanalytic set is subanalytic.

- The projection of a relatively compact subanalytic set is subanalytic.

- A finite union of subanalytic sets is subanalytic.

- A finite intersection of subanalytic sets is subanalytic.
\end{prop}
\begin{prop}\label{criterion}{\bf (\cite{BM} Proposition 3.13)}
Let $X$ be a real analytic manifold, and $M$ be a subset of $X$. Then the following conditions are
equivalent\,:

1) $M$ is subanalytic.

2) Every point of $X$ has a neighbourhood $U$ such that
$$M\cap U=\bigcup_{i=1}^p(f_{i1}(A_{i1})\setminus f_{i2}(A_{i2}))$$
where, for $i=1,\dots,p$ and $j=1,2$, $A_{ij}$ is a closed analytic subset of a real analytic manifold
$N_{ij}$, $f_{ij}\,:N_{ij}\rightarrow U$ is real analytic, and $f_{ij}\mid_{A_{ij}}\,:A_{ij}\rightarrow U$
is proper.

3) Every point of $X$ has a neighbourhood $U$ such that $M\cap U$ belongs to the class of subsets of
$U$ obtained using finite intersection, finite union and complement, from the family of closed subsets
of $U$ of the form $f(A)$ where $A$ is a closed analytic subset of a real analytic manifold $N$,
$f\,:N\rightarrow U$ is real analytic, and $f\mid_{A}$ is proper.
\end{prop}
As a direct corollary of this proposition, we have\,:
\begin{cor}\label{corodiff}
Let $X$ and $Y$ be two analytic manifolds and $\varphi\,:X\rightarrow Y$ be an analytic diffeomorphism.
Then, if $M\subset X$ is subanalytic, the set $\varphi(M)\subset Y$ is subanalytic.
\end{cor}
\begin{defn}
\rm
Let $X$ and $Y$ be two analytic manifolds, we say that a map $\varphi\,:X\rightarrow Y$ is {\it subanalytic}
if its graph is a subanalytic subset of $X\times Y$.
\end{defn}

As a projection of a relatively compact subanalytic set is subanalytic, we have directly\,:

\begin{prop}
The image of a relatively compact subanalytic set through a subanalytic map is subanalytic.
\end{prop}

Let $M\subset\P^{n+1}$ be a subanalytic set, we denote by $M_n$ the maximal subset of $M$
having the property to be a hypersurface in $\P^{n+1}$. 
\begin{defn}
\rm We denote by {\tenreg Sa} the set of subanalytic subsets of $\P^{n+1}$. We define the map
$G^{\Aff}$ to be\,:
$$\matrix{G^{\Aff}&\,:&\hbox{\tenreg Sa}&\longrightarrow&\hbox{\tenreg
P}\,(\P^{n+1})\cr&&M&\longmapsto&\overline{G^{\Aff}_{M_n}(M_n)}\cr}$$ The image of
$M\in\hbox{\tenreg Sa}\ $ through $G^{\Aff}$ is called the dual of $M$.                  
\end{defn} 

We can directly generalize the map $G^{\Aff}_M$ to subanalytic sets of $\P^{n+1}$.
This generalization gives more flexibility in some proofs.

\begin{defn}
\rm
Let $M\subset\P^{n+1}$ be a subanalytic set, we call {\it affine Gauss map of $M$} the map
$G^{\Aff}_M$ defined by\,:
$$\matrix{G^{\Aff}_M&\,:&M&\longrightarrow&{\hbox{\tenreg
P\,}}(\P^{n+1})\hfill\cr&&m&\longmapsto&\displaystyle\{\lim_k(L_{y_k}M_n)^{\perp}\mid
\{y_k\}_{k\in\N}\subset M_n,\cr&&&&\displaystyle\lim_k y_k=m\ {\rm and}\lim_kL_{y_k}M_n\
\rm{exists}\}\cr}$$
\end{defn}

As a subanalytic set admits a Whitney stratification, we have that the map 
$G^{\Aff}_M$ is well defined and that the following lemma is true. We remark that a point of $M$
could be sent to a family of affine planes.

\begin{lem}
Let $M$ be a subanalytic subset of $\P^{n+1}$, then we have $G^{\Aff}(M)=\bigcup_{m\in M}G^{\Aff}_M(m)$.
\end{lem}

Let us make the link between the definition of $G^{\Aff}_M$ for $M$ a hypersurface
of $\R^{n+1}$ and for $M$ a hypersurface of $\P^{n+1}$.
For all $1\leq i\leq n+2$ we denote by $U_i$ the open subset $\{[u_1,\dots,u_{n+2}]\mid
u_i\not=0\}\subset\P^{n+1}$ and we define a local chart on $U_i$\,:
$$\matrix{\xi_i&\,:&U_i&\longrightarrow&\R^{n+1}\cr&&[u_1,\dots,u_{n+2}]&\longmapsto&
[{u_1\over u_i},\dots,{u_{i-1}\over u_i},{u_{i+1}\over u_i},\dots,{u_{n+2}\over u_i}]\cr}$$
Let $M\subset\P^{n+1}$ be a hypersurface, then for all $i$ we have\,:
$$\matrix{G^{\Aff}_M\!\mid_{M\cap U_i}&\,:&M\cap
U_i&\rightarrow&\P^{n+1}\cr &&m&\mapsto&(L_mM)^{\perp}=[N(m),<N(m),\xi_i(m)>]\cr}$$ where
$N(m)$ is a normal to $\xi_i(M\cap U_i)$ at $\xi_i(m)$. So the definition of $G^{\Aff}_M$
for a hypersurface in $\P^{n+1}$ is the natural generalization of the map $G^{\Aff}_M$ for
$M$ a hypersurface in $\R^{n+1}$.

Now we begin with the first step to see that $G^{\Aff}$ is an involution on good subanalytic sets.
\begin{thm}{\label{premierthm}}
Let $M$ be a subanalytic subset of $\P^{n+1}$, then $G^{\Aff}(M)$ is subanalytic.
\end{thm}

To give the proof of theorem \ref{premierthm} we have to make some preparation work and we
will use the following auxiliary map\,:
$$\matrix{\alpha&\,:&\R^{n+1}\times\P^{n+1}&\rightarrow&\P^{n+1}\cr&&(m,[v])&\mapsto&
[v,<m,v>]}$$
\begin{lem}
The map $\alpha$ is subanalytic.
\end{lem}
\begin{pf}
As the graph of $\alpha$ is closed, it is sufficient to check the criterion at every point of the
graph. We consider the following map\,:
$$\matrix{\beta&\,:&\R^{n+1}\times\P^n&\rightarrow&\R^{n+1}\times\P^n\times\P^{n+1}\cr
&&(q,[w])&\longmapsto&(q,[w],[w,<q,w>])\cr}$$ This is a proper analytic map from an analytic
set, and as this map is surjective on the graph of $\alpha$ the criterion is
checked.               
\findemo\end{pf}

Now we have to use the tangent mapping theorem we can find in \cite{Loj2}\,:
\begin{thm}
{\bf (\cite{Loj2} p. 1590)} Let $M\subset\R^{n+1}$ be a subanalytic hypersurface, 
then the Gauss map of $M$ is subanalytic.
\end{thm}

\begin{prop}
Let $M\subset\P^{n+1}$ be a subanalytic hypersurface, then the map $G^{\Aff}_M$ is subanalytic.
\end{prop}
\begin{pf}
We have to prove the subanalycity in $\P^{n+1}\times\P^{n+1}$ of

$\{(m,(L_mM)^\perp)\in\P^{n+1}\times\P^{n+1}\mid m\in M\}$ .

Let $m\in M$ and $i\in\{1,\dots,n+2\}$ such that $m\in U_i$. We know that $M\cap U_i$ is subanalytic,
so using the corollary \ref{corodiff} we have that $\xi_i(M\cap U_i)$ is a subanalytic hypersurface of
$\R^{n+1}$. Now, using the tangent mapping theorem, we have that the set
$\{(\xi_i(m),[N(m)])\mid m\in M\cap U_i\}$, where $N(m)$ is a normal vector to
$\xi_i(M\cap U_i)$ at $\xi_i(m)$, is subanalytic in $\R^{n+1}\times\P^n$. As the map $\alpha$ is
subanalytic, we have the subanalycity of the set\,:
$$\{(\xi_i(m),[N(m)],[N(m),<N(m),\xi(m)>])\mid m\in M\cap U_i\}=$$
$$\{(\xi_i(m),[N(m)],(L_mM)^\perp)\mid m\in M\cap U_i\}$$
Now projecting this set and using the analytic diffeomorphism given by $\xi_i^{-1}\times\hbox{Id}$, we
have that $\{(m,(L_mM)^{\perp})\in (M\cap U_i)\times\P^{n+1}\}$ is subanalytic in $U_i\times\P^{n+1}$.
We have checked the subanalycity criterion in a neighbourhood of each point, so we have the conclusion.
\findemo\end{pf}

\begin{cor}
Let $M\subset\P^{n+1}$ be a subanalytic hypersurface, then $G^{\Aff}_M(M)$ is
subanalytic.
\end{cor}

\begin{pf}
The set $G^{\Aff}_M(M)$ is the image of a relatively compact subanalytic set by a subanalytic map, so
$G^{\Aff}_M(M)$ is subanalytic in $\P^{n+1}$.
\findemo\end{pf}

\begin{pf} (Theorem \ref{premierthm})
Let $M\subset\P^{n+1}$ be a subanalytic set, we have that $M_n$ is a subanalytic hypersurface of
$\P^{n+1}$. Then $G^{\Aff}_{M_n}(M_n)$ is subanalytic and as the closure of a subanalytic set is 
subanalytic, we have the subanalycity of $G^{\Aff}(M)=\overline{G^{\Aff}_{M_n}(M_n)}$.
\findemo\end{pf}
\section{Restriction of $G^{\Aff}$ to admissible sets}
We want a duality theorem on $G^{\Aff}$, and it is clear there exists subanalytic sets (for
instance a subanalytic set $M$ with $M_n=\emptyset$) such that $G^{\Aff}(G^{\Aff}(M))\not=M$. So we
have to restrict $G^{\Aff}$ to a smaller class of subsets of $\P^{n+1}$ than the subanalytic class.
\begin{defn}
\rm
Let $M\subset\P^{n+1}$ be a hypersurface and $m\in M$. We denote by $r(m)$ the rank of the map
$G^{\Aff}_M$ at $m$ and we define the subset $M^{[k]}$ of $M$ by\,:
$$M^{[k]}=\{m\in M\mid r(m)=k\}$$
\end{defn}
\begin{defn}
\rm
We say that $M\subset\P^{n+1}$ is {\it admissible} if we have\,:

1) $M$ is subanalytic.

2) $(M_n)^{[n]}$ is dense in $M$.

3) $M$ is closed.

We denote by {\tenreg Ad} the set of admissible subsets of $\P^{n+1}$.
\end{defn}
In this section we denote  by $U$ the open $U_1\subset\P^{n+1}$ and $\xi$ the map $\xi_1$. We remark
that if $M\in${\tenreg Ad}, then $(M_n\cap U)^{[n]}$ is dense in $M$, so it is sufficient to work in
$U$ for the following proofs.

Now we are going to prove that the image through $G^{\Aff}$ of an admissible set is admissible.

\begin{lem}
Let $M\in${\tenreg Ad}, then $(M_n\cap U)^{[n]}$ is subanalytic.
\end{lem}
\begin{pf}
We know that $M_n$ and $U$ are subanalytic, so $M_n\cap U$ is subanalytic too. We can apply 
proposition \ref{criterion} and for all $m\in M_n\cap U$ there exists $V$ a neighbourhood of $m$ and all
the things we need to write\,:
$$(M_n\cap U)\cap V=\bigcup_{i=1}^p(f_{i1}(A_{i1}\setminus f_{i2}(A_{i2}))$$
We now consider the maps $\delta_i\,:A_{i1}\rightarrow\R$ which send a point $a\in A_{i1}$ to the
determinant of the derivative of the Gauss map of $\xi(M_n\cap U)$ at $\xi(f_{i1}(a))$. We remark that
$\delta_i^{-1}(0)$ is an analytic closed subset of $N_{i1}$, denoted by $B_{i1}$, and
that$f_{i1}\mid_{B_{i1}}$ is proper. So we have\,:
$$(M_n\cap U)^{[n]}\cap V=\bigcup_{i=1}^p(f_{i1}(A_{i1}\setminus (f_{i1}(B_{i1})\cup f_{i2}(A_{i2})))$$
and the conclusion follows by proposition \ref{criterion}.
\findemo\end{pf}
 
\begin{lem}\label{lemlocal}
Let $M\subset U\subset\P^{n+1}$ be a subanalytic hypersurface which verifies $M=M_n=M^{[n]}$ and
$G^{\Aff}_M(M)=(G^{\Aff}_M(M))_n\subset U$. Then $G^{\Aff}_M$ is an analytic diffeomorphism on its
image whose inverse is $G^{\Aff}_{G^{\Aff}(M)}$, moreover $(G^{\Aff}_M(M))^{[n]}$ is dense in
$G^{\Aff}_M(M)$.
\end{lem}
\begin{pf}
As  the image of $M\subset U$ through $G^{\Aff}_M$ is contained in $U$, we can use the formula we have
in the case of hypersurfaces in $\R^{n+1}$ dividing by the non zero factor $<N(v),f(v)>$. We obtain\,:
$$\matrix{G^{\Aff}_M&\,:&M&\longrightarrow&G^{\Aff}_M(M)\cr&&&&\cr
&&m&\longmapsto&\big[{N(v)\over<N(v),f(v)>},-1\big]\cr}$$ where $f\,:V\rightarrow\R^{n+1}$
is a local parametrization of $\xi(M)$ at $\xi(m)$ and
$N(v)$ is a normal to $\xi(M)$ at $\xi(f(v))$. As $G^{\Aff}_M(M)=(G^{\Aff}_M(M))_n\subset
U$, we have that
${N(v)\over<N(v),f(v)>}$ is a local parametrization of
$\xi(G^{\Aff}_M(M))$ at $\xi(G^{\Aff}_M(f(v)))$. Let us prove that the matrix 
$$\pmatrix{{N(v)\over<N(v),f(v)>}\cr
\dd{v_i}({N(v)\over<N(v),f(v)>})\cr}$$
is invertible.

We suppose there exist $\lambda,\lambda_1,\dots,\lambda_n$ such that\,:
$$\lambda{N\over<N,f>}=\sum_{i=1}^n\lambda_i
\dd{v_i}({N\over<N,f>})\Longrightarrow$$
$$\lambda{N\over<N,f>}=\sum_{i=1}^n\lambda_i
({\dd{v_i}N<N,f>-N<\dd{v_i}N,f>\over<N,f>^2})\Longrightarrow$$
$$(\lambda+\sum_{i=1}^n{\lambda_i<\dd{v_i}N,f>\over<N,f>})N=\sum_{i=1}^n\lambda_i\dd{v_i}N$$
As $M=M^{[n]}$, the family
$\{N(v),\dd{v_1}N(v),\dots,\dd{v_n}N(v)\}$ is linearly independent for all $v$. So the $\lambda_i$
are all equal to $0$, it follows that $\lambda=0$. We have proved that the line vectors of the
matrix are linearly independent, so this matrix is invertible.

We can see that\,:
$$\matrix{G^{\Aff}_{G^{\Aff}_M(M)}&\,:&V&\longrightarrow&\P^{n+1}\cr
&&v&\longmapsto&\bigg[\pmatrix{{N(v)\over<N(v),f(v)>}\cr
\dd{v_i}({N(v)\over<N(v),f(v)>})\cr}^{-1}\pmatrix{1\cr
0_n\cr},1\bigg]\cr}$$
As the matrix of the equation\,:
$$\pmatrix{{N(v)\over<N(v),f(v)>}\cr\dd{v_i}({N(v)\over<N(v),f(v)>})\cr}X=\pmatrix{1\cr
0_n\cr}$$ is invertible, this equation has, for each $v$, a unique solution. It is easy to
see that $f(v)$ is a solution of this equation. It follows that\,:
$$G^{\Aff}_{G^{\Aff}_M(M)}(G^{\Aff}_M(m))=[f(v),1]=m$$
So we have $G^{\Aff}_{G^{\Aff}_M(M)}(G^{\Aff}_M(M))=M$. As $G^{\Aff}_M$ is a diffeomorphism on its
image, $G^{\Aff}_{G^{\Aff}_M(M)}$ is the inverse of $G^{\Aff}_M$. For the last point, if the set
$G^{\Aff}_M(M)^{[n]}$ is not dense in
$G^{\Aff}_M(M)$, then the map $G^{\Aff}_{G^{\Aff}_M(M)}$ is not bijective, it is due to the
fact that an open set of $G^{\Aff}_M(M)$ is sent to a part of strictly lower dimension than the
dimension of $M$.
\findemo\end{pf}
\begin{lem}\label{lemdensite}
Let $M\subset\P^{n+1}$ be a subanalytic set and $V\subset M$ be an open dense set in $M$, then
$G^{\Aff}(V)=G^{\Aff}(M)$.
\end{lem}
\begin{pf}
First of all, we prove this lemma for $M$ a hypersurface. It is clear that
$G^{\Aff}_V=G^{\Aff}_M\!\mid\!_V$, as $V$ is dense in $M$ and $G^{\Aff}_M$ is continuous, we have
$G^{\Aff}(M)=\overline{G^{\Aff}_M(M)}=\overline{G^{\Aff}_V(V)}=G^{\Aff}(V)$.

For the case of $M$ subanalytic, as $V$ is open dense in $M$, we have that $V_n$ is open dense in
$M_n$. We obtain\,:
$$G^{\Aff}(M)=\overline{G^{\Aff}_{M_n}(M_n)}=\overline{G^{\Aff}_{V_n}}=G^{\Aff}(V)$$.\findemo\end{pf}
\begin{prop}
Let $M\in${\tenreg Ad}, then $G^{\Aff}(M)\in${\tenreg Ad}.
\end{prop}
\begin{pf}
By definition of $G^{\Aff}$ we have that $G^{\Aff}(M)$ is closed and theorem \ref{premierthm} gives us
that $G^{\Aff}(M)$ is subanalytic. So the only fact it remains to prove is the density of
$(G^{\Aff}(M)_n)^{[n]}$ in $G^{\Aff}(M)$.

By definition, we have $M=\overline{(M_n\cap U)^{[n]}}$, so by lemma \ref{lemdensite}
we obtain $G^{\Aff}(M)=\overline{G^{\Aff}_M((M_n\cap U)^{[n]})}$.
We know that $(M_n\cap U)^{[n]}$ is subanalytic, that implies the subanalycity of $G^{\Aff}_M((M_n\cap
U)^{[n]})$. So the map $G^{\Aff}_M\!\mid\!_{(M_n\cap U)^{[n]}}=G^{\Aff}_{(M_n\cap U)^{[n]}}$ is
an immersion of $(M_n\cap U)^{[n]}$ in $\P^{n+1}$ whose image is subanalytic. Let us consider a locally
finite stratification of $G^{\Aff}_M((M_n\cap U)^{[n]})$ with analytic manifolds. Let $p\in
G^{\Aff}_M((M_n\cap U)^{[n]})$, as $G^{\Aff}_{(M_n\cap U)^{[n]}}$ is an immersion, each neighbourhood
of $p$ intersects a submanifold of dimension $n$, as the stratification is locally finite, $p$ is
adherent to a $n$-dimensional strata. So we have that $(G^{\Aff}((M_n\cap U)^{[n]}))_n$ is dense in
$G^{\Aff}_M((M_n\cap U)^{[n]})$. As $(G^{\Aff}((M_n\cap U)^{[n]}))_n\subset G^{\Aff}(M)_n$ and
$G^{\Aff}(M)=\overline{G^{\Aff}_M((M_n\cap U)^{[n]})}$, we have that $G^{\Aff}(M)_n$ is dense in
$G^{\Aff}$.

To see that $(G^{\Aff}(M)_n)^{[n]}$ is dense in $G^{\Aff}(M)$, let us consider the following subset of
$M$\,:
$$M'=\big(G^{\Aff}_{(M_n\cap U)^{[n]}}\big)^{-1}\big((G^{\Aff}_M((M_n\cap U)^{[n]}))_n\big)$$
The set $M'$ is the subset of $(M_n\cap U)^{[n]}$ containing the points whose tangent hyperplane is not
bitangent.

The regular part of maximal dimension of a subanalytic set $M$ is open in $M$, so we have that
$G^{\Aff}_M((M_n\cap U)^{[n]})_n$ is open in $G^{\Aff}_M(M)$. Moreover we know that 
$G^{\Aff}((M_n\cap U)^{[n]})_n$ is dense in $G^{\Aff}_M(M)$. As $G^{\Aff}_M$ restricted to $(M_n\cap
U)^{[n]}$ is locally a diffeomorphism, the set $M'$ is open dense in $(M_n\cap U)^{[n]}$. We have that
$M'$ is a dense submanifold of $(M_n\cap U)^{[n]}$ and $\overline{(M_n\cap U)^{[n]}}=M$, so we have
$\overline{M'}=M$.

The set $M'$ is the inverse image of a relatively compact subanalytic set through a subanalytic
map, so $M'$ is subanalytic. We can now consider the map $G^{\Aff}_{M'}=G^{\Aff}_M\!\mid\!_{M'}$, we can
apply the lemma \ref{lemlocal} to $M'$, so we have that $G^{\Aff}_{M'}$ is an analytic diffeomorphism on
its image and that $G^{\Aff}_{M'}(M')^{[n]}$ is dense in $G^{\Aff}_{M'}(M')$.

As $G^{\Aff}_{M'}(M')$ is dense in $G^{\Aff}(M)$, we have the density of $G^{\Aff}_{M'}(M')^{[n]}$ in
$G^{\Aff}_M(M)$, finally $G^{\Aff}_{M'}(M')^{[n]}\subset(G^{\Aff}_M(M)_n\cap U)^{[n]}$ and we have
$G^{\Aff}(M)\in${\tenreg Ad}.
\findemo\end{pf}
\begin{thm}\label{bigthm}
The map $G^{\Aff}$ restricted to {\tenreg Ad} is an involution.
\end{thm}
\begin{pf}
We have seen that $G^{\Aff}(M)=\overline{G^{\Aff}_M(M')}$ and $\overline{M'}=M$. So we have
$\overline{G^{\Aff}_{G^{\Aff}_{M'}(M')}(G^{\Aff}_{M'}(M'))}=G^{\Aff}\circ G^{\Aff}(M)$ and by lemma
\ref{lemlocal} we have $G^{\Aff}\circ G^{\Aff}(M)=\overline{M'}=M$.
\findemo\end{pf}

\section{Intersection of admissible sets with linear subspaces}

In this part, we are interested in the study of intersections between admissible sets of $\P^{n+1}$ and
linear subspaces of codimension 1 (in other words hyperplanes). We will prove that if $M\subset\P^{n+1}$
is admissible, then there exists a dense set $W$ of hyperplanes with $M\cap\sigma$ admissible for all
$\sigma\in W$.

Let $M\subset\P^{n+1}$ be admissible and $\sigma\subset\P^{n+1}$ be a hyperplane, then $M\cap\sigma$ is
closed and subanalytic. So the point to study is for which $\sigma$ we have the density of
$((M\cap\sigma)_{n-1})^{[n-1]}$ in $M\cap\sigma$.

First we need the following lemma whose proof is straight forward\,:

\begin{lem}
Let $f\,:\R^n\rightarrow\R^k$ be a differentiable map with $f(x)\not=0$ for all $x\in\R^n$. We have\,:
$$\hbox{Ker}(D({f\over\Vert f\Vert})(x))=\hbox{Ker}Df(x)\oplus\hbox{vect}\{Df(x)^{-1}(f(x))\}$$
where $Df$ denotes the derivative, $\hbox{vect}\{\,\}$ denotes the vector space generated by $\{\,\}$ and
$T_{f(x)}\R^k$ is identified with $\R^k$.
\end{lem}

Let us do some work on hypersurfaces of $\R^{n+1}$, so let $M\subset\R^{n+1}$ be such a
hypersurface, we denote by
$M^{[n]}$ the set of points in
$M$ where the Gauss map is of rank $n$. Given a hypersurface $M\subset\R^{n+1}$, we will prove there
exists a dense set
$D\subset\Aff(n,1)$ such that $(M\cap\tau)^{[n-1]}=M\cap\tau$ for all $\tau\in D$, which is in fact more
than we need. Then we will obtain the result on admissible sets as a corollary.

First we will determine in which case we loose the property that the Gauss map is of maximal rank by
intersecting with a hyperplane.

\begin{prop}\label{prolocal}
Let $M\subset\R^{n+1}$ be a hypersurface with $M^{[n]}=M$ and $p\in M$. We consider $\sigma\in G(n,1)$
which intersects $T_pM$ transversally and we denote by $\tau\in\Aff(n,1)$ the hyperplane parallel to
$\sigma$ through $p$.

Then the Gauss map of $\tau\cap M\subset\tau$ at $p$ is not of maximal rank if and only if the second
fundamental form of $M$ at $p$ restricted to $T_p(M\cap\tau)$ is degenerated.
\end{prop}

\begin{pf}
As we work locally, we can suppose that the Gauss map has values in $S^n$, the $n$-dimensional sphere.
Using an isometry of 
$\R^{n+1}$, we can suppose that
$\tau=\R^n\times\{0\}$ and $p=0$. As $\sigma$ is transverse to $T_pM$, there exists a neighbourhood $U$ of
$p$ in $\tau$ and a neighbourhood $V$ of $p$ in $\R^{n+1}$ with
$$\matrix{\varphi&\,:&U\subset\tau=\R^n\times\{0\}&\longrightarrow&V\subset\R^{n+1}\cr
                 &&(x_1,\dots,x_n)&\longmapsto&(x_1,\dots,x_n,h(x_1,\dots,x_n))\cr}$$
a local parametrization of $M$ at $p$ and $h$ a smooth function with ${\rm grad}\
h\not=0$ on
$U$. The parametrization of
$M$ we use allows us to see $M$ as the graph of the function $h$. The Gauss map in local
coordinates is given by\,:
$$\matrix{G_M\circ\varphi&\,:&U&\longrightarrow&S^n\cr& &x&\longmapsto &\displaystyle{({\rm
grad}\ h(x),-1)\over\Vert({\rm grad}\  h(x),-1)\Vert}\cr}$$ We consider, as a chart for
$S^n$, the central projection of the open hemisphere which contains
$(0_{n},-1)$ on the tangent plane of $S^n$ at
$(0_{n},-1)$. In this chart, denoted by $\psi$, we have $\psi\circ
G_M\circ\varphi(x)={\rm grad}\ h(x)$. As $M\cap\tau=h^{-1}(0)$, we have that ${\rm grad}\ h(s)$ is
orthogonal to $T_s(M\cap\tau)$ for all $s\in M\cap\tau\cap V$. So we have\,:
$$\matrix{G_{M\cap\tau}&\,:&M\cap\tau\cap V&\longrightarrow &S^{n-1}\cr & &s&\longmapsto
&\displaystyle {{\rm grad}h(s)\over\Vert {\rm grad}h(s)\Vert}\cr}$$
As $M=M^{[n]}$, the map ${\rm grad}\ h$ is smooth of maximal rank. From the preceding lemma we have\,:
$${\rm Ker}DG_{M\cap\tau}(p)={\rm vect}\{(D{\rm grad}\ h(p))^{-1}({\rm grad}\ h(p))\}\cap T_p(M\cap\tau)$$
Finally we have the following\,:
$${\rm Ker}DG_{M\cap\tau}(p)\not=\{0\}\Longleftrightarrow $$
$$\exists v\in T_p(M\cap\tau)\
{\rm s.\ t.}\ DG_M(p)(v)=\lambda{\rm grad}h(p)\Longleftrightarrow$$
$$\exists v\in T_p(M\cap\tau)\ {\rm s.\ t.}\ DG_M(p)(v)\perp
T_p(M\cap\tau)\Longleftrightarrow$$
$$\exists v\in T_p(M\cap\tau)\ {\rm s.\ t.}\
<DG_M(p)v,w>=0\ \forall w\in T_p(M\cap\tau)\Longleftrightarrow$$
$$\hbox{I\hskip-2pt I}_p\mid_{T_p(M\cap\tau)}\ {\rm is\ degenerate}$$
where $\hbox{I\hskip-2pt I}_p$ denotes the second fundamental form of $M$ at $p$.
\findemo\end{pf}

We will use the following transversality lemma\,:

\begin{lem}{\bf (\cite{Gib} p. 49)}\label{Gibson} Let $N$, $S$ and $P$ be manifolds and consider
$F\,:N\times S\rightarrow P$ a smooth family of smooth mappings transverse to smooth submanifolds $Q_1,
\dots, Q_t$ of $P$. Then there is a dense set of parameters $s\in S$ for which $F\mid_{N\times\{s\}}$ is
transverse to all of $Q_1,
\dots, Q_t$.
\end{lem}

\begin{defn}
\rm 
Let $M\subset\R^{n+1}$ be a hypersurface with $M=M^{[n]}$. {\it The isotropic cone of $M$}
is the set\,:
$$C_M=\{(m,x)\mid m\in M, x\in T_mM\ \hbox{with}\ \hbox{I\hskip-2pt I}_m(x,x)=0\ \hbox{and}\
x\not=0\}$$$$\subset TM\subset\R^{n+1}\times\R^{n+1}$$
where $\hbox{I\hskip-2pt I}_m$ is the second fundamental form of $M$ at $m$.
\end{defn}

We remark that $C_M$ is a submanifold of $\R^{n+1}\times\R^{n+1}$. If the second fundamental form of $M$
is everywhere definite (positive or negative), we have $C_M=\emptyset$. In the other case, $C_M$ has
codimension $3$. The notion of isotropic cone of $M$ will be important.

\begin{lem}
Let $M\subset\R^{n+1}$ be a hypersurface and $\tau\in\Aff(n,1)$, we denote by $T\tau\subset T\R^{n+1}$ the
tangent bundle of $\tau$. If $T\tau$ is transverse to $C_M$ in $T\R^{n+1}$ and if $\tau$ is
tranverse to $M$ in
$\R^{n+1}$, then $\tau\cap M$ is a submanifold of $\tau$ and $(\tau\cap M)^{[n-1]}=\tau\cap M$.
\end{lem}

\begin{pf}
Because of the transversality of $M$ and $\tau$, we have that $\tau\cap M$ is a submanifold. For the second
affirmation we remark that if $\varphi$ is a non-degenerate quadratic form of $\R^n$ and if $W\subset\R^n$
is a subspace of codimension $1$, then $\varphi\mid_W$ is degenerate if and only if $W$ is tangent to
the isotropic cone of $\varphi$. As $\tau$ is transverse to $M$, the transversality of $T\tau$ and $C_M$
implies that at each point $p\in \tau\cap M$ the plane $T_p\tau\cap T_pM$ is not tangent to the isotropic
cone of the second fundamental form of $M$ at $p$. Proposition \ref{prolocal} allows us to conclude.
\findemo\end{pf}

Now we will construct an auxiliary map to apply lemma \ref{Gibson}. It is easier to work in
$\R^{n+1}\times G(n,1)$ than in $\Aff(n,1)$. So, we will consider, after using this auxiliary map, the
canonical submersion of $\R^{n+1}\times G(n,1)$ on $\Aff(n,1)$.

Let us denote by $\gamma$ the total space of the canonical bundle on $G(n,1)$, we define\,:
$$\matrix{g&\,:&\R^{n+1}\times(\gamma\oplus\gamma)&\longrightarrow&T\R^{n+1}\cr
           &&(z,(\tau,x,y))&\longmapsto&(x+z,y)\cr}$$

\begin{prop}
Let $M\subset\R^{n+1}$ be a hypersurface, then there exists a dense set $D\subset\Aff(n,1)$ such that
$\tau$ and $M$ are transverse for all $\tau\in D$.
\end{prop}

\begin{pf}
Let $\sigma=\R^n\times\{0\}\in G(n,1)$, we will use the following charts for $G(n,1)$.
We define $U_\sigma=\{\tau\in G(n,1)\mid \tau\cap\sigma^\perp=\{0\}\}$ a neighbourhood of $\sigma$,
and we consider the basis of $\sigma$ (seen as a vector subspace of $\R^{n+1}$) given by
$B=(I_n\vert 0)\in M_{n\times(n+1)}(\R)$ and the basis of $\sigma^{\perp}$ given by $B^{\perp}=(0\dots 0\
1)\in M_{1\times(n+1)}(\R)$. The map defined by\,:
$$\matrix{\varphi_{B,B^{\perp}}&\,:& M_{n\times 1}(\R)&\longrightarrow&U_\sigma\cr
                               &&A&\longmapsto&q(B+AB^{\perp})\cr}$$
where $q$ is the canonical projection of the set of rank $n$ matrices in $M_{n\times(n+1)}(\R)$ on
$G(n,1)$, gives a local chart of $G(n,1)$. 

Let us consider the map, denoted by $f$, which is the local description of $g$\,:
$$\matrix{\R^n\times\R^n\times\R^{n+1}\times U_\sigma&\longrightarrow&\R^{n+1}\times\R^{n+1}\cr
           (x,y,z,\varphi_{B,B^{\perp}}(A))&\longmapsto&((B+AB^{\perp})^tx+z,(B+AB^{\perp})^ty)\cr}$$
The derivative of $f$ at the point $(x,y,z,A)$ is equal to\,:
$$\pmatrix{(B+AB^{\perp})^t& 0_{(n+1)\times n}& I_{n+1} &\pmatrix{0_{n\times
n}\cr\overline{\matrix{x_1&\dots&x_n\cr}}}\cr
0_{(n+1)\times n}&(B+AB^{\perp})^t&0_{(n+1)\times(n+1)}&\pmatrix{0_{n\times
n}\cr\overline{\matrix{y_1&\dots&y_n\cr}}}\cr}$$
We will proceed in two steps to study the rank of $Df$. The $n+1$
first lines are linearly independent and as $B+AB^{\perp}$ has rank
$n$ we have that
$Df$ has rank greater than $2n+1$. Let us denote by $A$ the matrix $(a_1\dots a_n)$, we remark that
$$\pmatrix{(B+AB^{\perp})^t&\pmatrix{0_{n\times n}\cr\overline{\matrix{y_1&\dots&y_n\cr}}}\cr}=$$$$
   \pmatrix{\pmatrix{I_n\cr\overline{\matrix{a_1&\dots&a_n\cr}}}&\pmatrix{0_{n\times
n}\cr\overline{\matrix{y_1&\dots&y_n\cr}}}\cr}$$
has rank $n+1$ if $y\not=0$ and rank $n$ if $y=0$. In brief we have\,:
$$\hbox{rang}Df(x,y,z,A)=\left\{\matrix{2n+2&\hbox{si}\ y\not=0\cr 2n+1&\hbox{si}\
y=0\cr}\right.$$ 
So $f\mid\!_{\R^n\times(\R^n\setminus\{0\})\times\R^{n+1}\times U_\sigma}$ is a submersion, in other words,
we have the transversality of $f\mid\!_{\R^n\times(\R^n\setminus\{0\})\times\R^{n+1}\times
U_\sigma}$ and $C_M$ in $\R^{n+1}\times\R^{n+1}$. But as
$f(\R^n\times\{0\}\times\R^{n+1}\times U_\sigma)\cap C_M=\emptyset$ we have the transversality of
$f$ and $C_M$.

Moreover, it is easy to see that $f$ is transverse to $(M\times\R^{n+1})$, so there exists a dense set
$W\subset\R^{n+1}\times U_\sigma$ such that $f_{z,\tau}\,:(x,y)\mapsto f(x,y,z,\tau)$ is transverse to
$C_M$ and $M\times\R^{n+1}$ for all $(z,\tau)\in W$. Or the image of $f_{z,\tau}$ is the affine plane
parallel to $\tau$ containing $z$, so we have a dense part $V\subset\Aff(n,1)$ such that for all
$\pi\in V$ we have $T\pi$ transverse to $C_M$ and $T\pi$ transverse to $M\times\R^{n+1}$.
\findemo\end{pf}

\begin{cor}
Let $M\subset\P^{n+1}$ be an admissible set. There exists a dense set $W$ of linear subspaces of
codimension $1$ in $\P^{n+1}$ with $M\cap\sigma$ admissible for all $\sigma\in W$.
\end{cor}

\begin{pf}
We denote by $U$ the open $U_{n+2}\subset\P^{n+1}$ and by $\xi$ the map
$\xi_{n+2}\,:U_{n+2}\rightarrow\R^{n+1}$. We can apply the preceding proposition to
$\xi(M_n^{[n]}\cap U)$. So there exists a dense set $V\subset\Aff(n,1)$ such that for all
$\pi\in V$ we have
$(\pi\cap\xi(M_n^{[n]}\cap U))^{[n-1]}=\pi\cap\xi(M_n^{[n]}\cap U)$. We remark that $\xi(\pi)$ is a linear
subspace of codimension $1$ in
$\P^{n+1}$ and that $\xi^{-1}(\pi\cap\xi(M\cap U))=\xi^{-1}(\pi)\cap M\cap U$. Then we obtain
$(\xi(\pi)\cap M\cap U)_n^{[n]}=\xi(\pi)\cap (M_n)^{[n]}\cap U$ and
$M\cap\xi(\pi)$ is admissible. The density of $V$ in $\Aff(n,1)$ implies the density of
$\{\xi(\pi)\mid\pi\in V\}$ in the set of linear subspaces of codimension $1$ in
$\P^{n+1}$; so we have the conclusion.
\findemo\end{pf}

\section{Duality and outlines of hypersurfaces}

In this section we will study the relations between outlines of a compact analytic hypersurface $M$ in
$\R^{n+1}$ and the dual of $M$. 

Let $\sigma\in\P^{n+1}$, we denote by 
$\P(\sigma^\perp)$ the linear subspace of $\P^{n+1}$ containing all the normal directions to $\sigma$.
As $\P(\sigma^\perp)$ is isomorphic to $\P^n$, we can apply
$G^{\Aff}$ to subanalytic subspaces of $\P(\sigma^\perp)$ and we can consider their image through $G^{\Aff}$
in $\P(\sigma^\perp)$.

We denote by $i\,:\R^{n+1}\rightarrow\P^{n+1}$ the injection given by $i(x)=[x,1]$. Let
$M\subset\R^{n+1}$ be a compact analytic hypersurface, then $i(M)$ is admissible and $G^{\Aff}(i(M))$ is
admissible too. We will construct the following diagram for $M$ a compact analytic
hypersurface in $\R^{n+1}$ and
$[v]\in\P^n$ a direction\,:
$$\matrix{M_{}&\hfl{G^{\Aff}\circ i}{}&G^{\Aff}(i(M))\cr
          \vfl{\pi_{[v]}}{}&&\vfl{}{\cap\P([v,0]^\perp)}\cr
          C^M_{[v]}&\hfll{\xi\circ G^{\Aff}}{}&G^{\Aff}(i(M))\cap\P([v,0]^\perp)\cr}$$
The outline $C^M_{[v]}$ is seen in the plane orthogonal to the direction $[v]$, so its image through $i$
is contained in $\P([v,0]^\perp)$. Then $G^{\Aff}(i(C^M_{[v]}))$ is the dual of
$i(C^M_{[v]})$ seen as a part of $\P([v,0]^\perp)$. There are two difficulties with the map
$\xi\circ G^{\Aff}$. The first one is that
$G^{\Aff}(i(C^M_{[v]}))$ is not necessarily equal (but it is included in) to
$G^{\Aff}(i(M))\cap\P([v,0]^\perp)$. The second one is due to the fact that an outline is not
necessarily admissible. In other words the involution property is not true for
$G^{\Aff}$ in the outline level, but this property is true for a dense set of directions.

We begin with some facts that are true for admissible sets in $\P^{n+1}$. 
Let $M$ be an admissible set in $\P^{n+1}$ and $\sigma\in\P^{n+1}$ be a direction. To obtain a link
between $G^{\Aff}(M)\cap\P{(\sigma^\perp)}$ and $G^{\Aff}(M\cap\P(\sigma^\perp))$, we need some auxiliary
map.

Let us consider the retraction of $\P^{n+1}\setminus\{\sigma\}$ on
$\P(\sigma^\perp)$ defined by\,:
$$\matrix{\rho_{\sigma}&\,:&\P^{n+1}\setminus\{\sigma\}&\longrightarrow &\P(\sigma^\perp)\cr
                         &&[x]&\longmapsto &[\pi_{\sigma}(x)]\cr}$$
This allows us to define a map, denoted by $\rho_\sigma$ too, from
$\hbox{\tenreg P}\,(\P^{n+1})$, the power set of $\P^{n+1}$, in $\hbox{\tenreg
P}\,(\P(\sigma^\perp))$\,:
$$\matrix{\rho_{\sigma}&\,:&\hbox{\tenreg P}\,(\P^{n+1})&\longrightarrow &\hbox{\tenreg
P}\,(\P(\sigma^\perp))\cr
      &&A&\longmapsto &\left\{\matrix{\cup_{a\in A}\rho_\sigma(a)&
                                       \hbox{si $\sigma\not\in A$}\cr
                                       \P(\sigma^\perp)&
                                       \hbox{si $\sigma\in A$}\cr}\right.\cr}$$
Now let us compare the maps $G^{\Aff}_M$ and $G^{\Aff}_{M\cap\P(\sigma^\perp)}$.

\begin{defn}\rm
Let $M$ be a subanalytic subset of $\P^{n+1}$ and $\sigma\in\P^{n+1}$, we denote by $M(\sigma)$ the set
$M\cap\P(\sigma^{\perp})$.
\end{defn}

\begin{lem}
Let $M$ be an admissible subset of $\P^{n+1}$ and $\sigma\in\P^{n+1}$, then for all
$m\in M(\sigma)$\,:
$$G^{\Aff}_{M(\sigma)}(m)\subset\rho_{\sigma}(G^{\Aff}_M(m))$$
Moreover we have the equality if $m$ is adherent to $M(\sigma)_{n-1}=(M(\sigma))_{n-1}$ and if
$\sigma\not\in G^{\Aff}_M(m)$. 
\end{lem}
\begin{pf}
The set $M$ is considered to be Whitney stratified.

First we prove the result for $m\in M(\sigma)_{n-1}$ with $\sigma\not\in G^{\Aff}_M(m)$.

Let $\{y_k\}_{k\in{\rm I\hskip -3pt N}}\subset M_n$ be a sequence contained in a $n$-dimensional connected
strata with $\lim_k y_k=m$ and $\lim_kL_{y_k}M_n$ exists. By definition of a Whitney stratification, we
have $\lim_k L_{y_k}M_n\supset L_mM(\sigma)_{n-1}$, so\,:
$$\lim_k(L_{y_k}M_n)^\perp=\big(\lim_kL_{y_k}M_n\big)^\perp\subset\big(L_mM(\sigma)_{n-1})^\perp$$ 
But $G^{\Aff}_{M(\sigma)}(m)$ is by definition equal to the orthogonal of
$L_mM(\sigma)_{n-1}$ in
$\P(\sigma^\perp)$, so it is equal to the element of $\P^{n+1}$ orthogonal to
$L_mM(\sigma)_{n-1}$ and to $\sigma$. In other words, $G^{\Aff}_{M(\sigma)}(m)$ is the
intersection of
$(L_mM(\sigma)_{n-1})^\perp$ with $\P(\sigma^\perp)$, so it is the image through $\rho_\sigma$ of
$\lim_k(L_{y_k}M_n)^\perp$. Finally
$\rho_{\sigma}(\lim_kG^{\Aff}_M(y_k))=G^{\Aff}_{M(\sigma)}(m)$,
so
$\rho_{\sigma}(G^{\Aff}_M(m))=G^{\Aff}_{M(\sigma)}(m)$.

Let us suppose now that $m$ is adherent to $(M(\sigma))_{n-1}$ and that 

$\sigma\not\in G^{\Aff}_M(m)$. Let
$\{y_k\}_{k\in{I\hskip -3pt N}}\subset (M(\sigma))_{n-1}$ be a sequence contained in a
$(n-1)$-dimensional connected strata with $\lim_k y_k=m$. We can suppose that
$\sigma\not\in G^{\Aff}_M(y_k)$ for all $k$.
We know, due to the first part of the proof, that
$\rho_{\sigma}(G^{\Aff}_M(y_k))=G^{\Aff}_{M(\sigma)}(y_k)$ for all $k$, 
so the equality is true when we take the limit at both sides. Using the definition of $G^{\Aff}_M$ for $M$
a subanalytic set we obtain\,:
$$\rho_{\sigma}(G^{\Aff}_M(m))=G^{\Aff}_{M(\sigma)}(m)$$
If $m\in M\cap P(\sigma^\perp)$ is not adherent to $(M(\sigma))_{n-1}$, then
$G^{\Aff}_{M(\sigma)}(m)=\emptyset$ and the result is clear.

If $m\in M\cap P(\sigma^\perp)$ is such that $\sigma\in G^{\Aff}_M(m)$, we have
$\rho_\sigma(G^{\Aff}_M(m))=\P(\sigma^\perp)$ and the result is clear.
\findemo\end{pf}

We remark that the hypothesis ``$M$ admissible'' is necessary, as a counter-example
for $M$ only subanalytic, we can take a closed plane curve in $\P^3$.

Now let us make the link between the dual of a compact analytic hypersurface of $\R^{n+1}$ and its
outline. Let $M$ be a compact analytic hypersurface of $\R^{n+1}$ and $[v]\in\P^{n}$, we recall
that $G^{\Aff}(i(C^M_{[v]}))$ is the dual of $i(C^M_{[v]})$ seen as a part of
$\P([v,0]^\perp)$. The set $G^{\Aff}(i(C^M_{[v]}))$ is contained in the set of tangent affine planes of $M$
which contain the direction $[v]$, for a dense set of directions we have the equality.
In other words we have\,:
$G^{\Aff}(i(C^M_{[v]}))\subset G^{\Aff}(i(M))\cap\P([v,0]^\perp)$.

\begin{prop}
Let $M$ be a compact analytic hypersurface of $\R^{n+1}$ and $[v]\in\P^{n}$, then\,:
$$i(C^M_{[v]})\supset G^{\Aff}\big(G^{\Aff}(i(M))\cap\P([v,0]^\perp)\big)$$
moreover if $i(C^M_{[v]})$ is admissible, then we have equality.
\end{prop}

\begin{pf}
We recall that for $[v]\in\P^n$ we denote by $\pi_{[v]}$ the orthogonal projection of $\R^{n+1}$ on
$v^{\perp}$. In our situation we have for all $[x]\in U$\,:
$$\pi_{[v]}(\xi([x]))=\xi\circ\rho_{[v,0]}(x)$$
For, $\rho_{[v,0]}([x_1,\dots,x_{n+2}])=[x_1-\lambda v_1,\dots,x_{n+1}-\lambda v_{n+1},x_{n+2}]$
with
$$\lambda={\sum_{i=1}^{n+1} x_iv_i\over\sum_{i=1}^{n+1} v_i^2}$$
that implies $\xi\circ\rho_{[v,0]}([x])=({x_1-\lambda v_1\over x_{n+2}},\dots,{x_{n+1}-\lambda v_{n+1}\over
x_{n+2}})=\pi_{[v]}(\xi([x]))$.

Let us prove the inclusion. Let $t\in G^{\Aff}(G^{\Aff}(i(M))\cap\P([v,0]^\perp))$ and consider
$Y=\{y\in G^{\Aff}(i(M))\cap\P([v,0]^\perp)\mid t\in
G^{\Aff}_{G^{\Aff}(i(M))\cap\P([v,0]^\perp)}(y)\}$ which is a non empty set. Let $X=\bigcup_{y\in
Y} G^{\Aff}_{G^{\Aff}(i(M))}(y)$, we have $X\subset i(M)$ and $L_xi(M)\subset\P([v,0]^{\perp})$ for all
$x\in X$. We obtain\,:
$$t\in\bigcup_{y\in Y}G^{\Aff}_{G^{\Aff}(i(M))\cap\P([v,0]^\perp)}(y)\subset\bigcup_{y\in Y}\rho_{[v,0]}G^{\Aff}_{G^{\Aff}(i(M))}(y)=$$
$$\bigcup_{x\in X}\rho_{[v,0]}(x)=i\big(\bigcup_{x\in X}\pi_{[v]}(\xi(x))\big)\subset i(C^M_{[v]})$$
Let us suppose that $i(C^M_{[v]})$ is admissible, we have\,:
$$G^{\Aff}(i(C^M_{[v]}))\subset G^{\Aff}(i(M))\cap\P([v,0]^\perp)$$
then, applying $G^{\Aff}$ to both sides, we have\,:
$$i(C^M_{[v]})=G^{\Aff}(G^{\Aff}(i(C^M_{[v]})))\subset G^{\Aff}(G^{\Aff}(i(M))\cap\P([v,0]^\perp))$$
which is the other inclusion.
\findemo\end{pf}

We remark that $i(C^M_{[v]})$ admissible and $C^M_{[v]}$ generic outline are not equivalent notions. 
But for a fixed $M$ there exists a dense (open dense for $M$ compact) set of directions $[v]$ such that
$i(C^M_{[v]})$ is admissible and $C^M_{[v]}$ is a generic outline.
In other words, for a dense set of directions $[v]$, we have\,:
$$G^{\Aff}(i(C^M_{[v]}))=G^{\Aff}(i(M))\cap\P([v,0]^\perp)$$
$$i(C^M_{[v]})=G^{\Aff}\big(G^{\Aff}(i(M))\cap\P([v,0]^\perp)\big)$$

\section{Duality and projective transformation}
 Now we can study the behaviour of the dual when we apply certain transformations on an admissible set.
This problem is quite difficult, but we can give a fully answer in the case of projective maps.
\begin{defn}
\rm
We say that the map $\alpha\,:\P^{n+1}\rightarrow\P^{n+1}$ is {\it projective } if there exists a linear
bijective map $A\,:\R^{n+2}\rightarrow\R^{n+2}$ with $\alpha([v])=[Av]$. We say that $A$ {\it represents} 
$\alpha$ or $\alpha$ is {\it represented} by $A$.
\end{defn}
\begin{prop}\label{propproj}
Let $M\subset\P^{n+1}$ be a subanalytic subset and 

$\alpha\,:\P^{n+1}\rightarrow\P^{n+1}$ be a projective map.
Then we have\,:
$$\alpha^*(G^{\Aff}(M))=G^{\Aff}(\alpha(M))$$
where $\alpha^*$ is represented by $(A^t)^{-1}$ if $\alpha$ is represented by $A$.
\end{prop}

\begin{pf}
First we prove the following for $M$ a hypersurface in $\P^{n+1}$\,:
$$\alpha^*(G^{\Aff}_M(m))=G^{\Aff}_{\alpha(m)}(\alpha(m))$$
We then have the equalities\,:
$$\alpha^*(G^{\Aff}_M(m))=\alpha^*(\{[v]\in\P^{n+1}\mid <v,w>=0\ \hbox{for all}\ [w]\in L_mM\})=$$
$$\{[v]\in\P^{n+1}\mid <A^tv,w>=0\ \hbox{for all}\ [w]\in L_mM\})=$$
$$\{[v]\in\P^{n+1}\mid <v,Aw>=0\ \forall [w]\in L_mM\})=$$
$$\{[v]\in\P^{n+1}\mid <v,w>=0\ \forall [w]\in \alpha(L_mM)\})=$$
$$\{[v]\in\P^{n+1}\mid <v,w>=0\ \forall [w]\in
L_{\alpha(m)}\alpha(M)\})=G^{\Aff}_{\alpha(M)}(\alpha(m))$$
Now in the case of a subanalytic subset $M$ of $\P^{n+1}$, we have\,:
$$G^{\Aff}(\alpha(M))=\overline{G^{\Aff}_{\alpha(M_n)}(\alpha(M_n))}=\overline{\alpha^*G^{\Aff}_{M_n}(M_n)}
=\alpha^*G^{\Aff}(M)$$
\findemo\end{pf}

\begin{cor}
Let $M,N\subset\R^{n+1}$ be two compact analytic hypersurfaces. If for all $[w]\in\P^n$ we have
that
$C^M_{[w]}$ and $C^N_{[w]}$ are homothetic of center $0$ in $w^{\perp}$, then $M$ and $N$ are homothetic of
center $0$ in $\R^{n+1}$.
\end{cor}

\begin{pf}
Denote by $W$ a dense set of directions such that the sets 

$\P([w,0]^{\perp})\cap G^{\Aff}(i(M))$ and
$\P([w,0]^{\perp})\cap G^{\Aff}(i(N))$ are admissible for all $[w]\in W$. As the sets
$\P([w,0]^{\perp})\cap G^{\Aff}(i(M))$ are admissible, we have that
$G^{\Aff}(i(C^M_{[w]}))=\P([w,0]^{\perp})\cap G^{\Aff}(i(M))$ for all $[w]\in W$. By hypothesis, for
all $[w]\in W$ there exists $\alpha_{[w]}\,:\P^n\rightarrow\P^n$ a projective map represented by
$\pmatrix{\lambda_{[w]}I_{n}&0\cr 0&1\cr}$ such that $\alpha_{[w]}\big(i(C^M_{[w]})\big)=i(C^N_{[w]})$.
So we have $\alpha^*_{[w]}(G^{\Aff}(i(C^M_{[w]})))=G^{\Aff}(i(C^N_{[w]}))$ and as
$\alpha_{[w]}$ and $\alpha_{[v]}$ agree on $\P([w,0]^{\perp})\cap\P([v,0]^{\perp})$ for all $[v],[w]\in W$,
we have that $\lambda_{[w]}=\lambda_{[v]}=\lambda$. Let us denote by $\alpha\,:\P^{n+1}\rightarrow\P^{n+1}$
the projective map represented by $\pmatrix{\lambda I_{n+1}&0\cr 0&1\cr}$, then $\alpha\big( 
i(C^M_{[w]})\big)=i(C^N_{[w]})$ for all $[w]\in W$. By proposition \ref{propproj} we have that
$\alpha^*(G^{\Aff}(i(C^M_{[w]})))=\alpha^*(\P([w,0]^{\perp})\cap G^{\Aff}(i(M)))=\P([w,0]^{\perp})\cap
G^{\Aff}(i(N))$.  But the density of $W$ gives
$G^{\Aff}(i(P))=\overline{\cup_{[w]\in W}\P([w,0]^{\perp})\cap G^{\Aff}(i(P))}$ for $P=M$ or $N$, so by
continuity of
$\alpha^*$ we have $\alpha^*(G^{\Aff}(i(M)))=G^{\Aff}(i(N))$. Once again by proposition \ref{propproj} we
have $\alpha(i(M))=i(N)$, in other words, $M$ and $N$ are homothetic of center $0$.
\findemo\end{pf}

This result is still true for compact hypersurfaces which are not necessarily analytic. We just need to suppose
that the hypersurface $M$ is subanalytic, which mean in particular that $M$ can contained peaces of $l$-planes
for $l\leq n$. First we can see, using the duality process, that
$M^{[n]}$ and $N^{[n]}$ are homothetic of center $0$, where $M^{[n]}$ denotes the set of points in $M$
where the Gauss map is of maximal rank. Then we use the reconstruction process given in \cite{Poi} to
obtain that $M$ and $N$ are homothetic of center $0$.

\vtab
The author would like to thank Oscar Burlet who was his advisor and Fran\c{c}ois Haab for useful
discussions.

\end{document}